\documentclass[a4paper,12pt,reqno]{amsart}
\usepackage{amssymb}\theoremstyle{plain}

\usepackage{amscd}
\usepackage{enumitem}

\newtheorem*{theorem*}{Theorem}

\newtheorem*{corollary*}{Corollary}

\newtheorem*{lemma*}{Lemma}

\newtheorem*{proposition*}{Proposition}

\newtheorem*{conjecture*}{Conjecture}
\theoremstyle{definition}
\newtheorem{definition}{Definition}
\newtheorem*{definition*}{Definition}
\theoremstyle{remark}
\newtheorem{remark}{Remark}
\newtheorem*{remark*}{Remark}
\newtheorem{example}{Example}
\newtheorem*{problem*}{Problem}
\newtheorem{problem}{Problem}

\begin{document}

\begin{center}
\title[QuantMat]{\bf {On non-Archimedean quantum mathematics and non-Archimedean (quantum) computation}\\
(spaces of non-Archimedean (quantum) mathematics)}
\end{center}

\maketitle

\begin{center}
{\bf Nikolaj M. Glazunov } \end{center}

\begin{center}
{\rm Glushkov Institute of Cybernetics NASU, Kiev, } \\

{\rm  Email:} {\it glanm@yahoo.com }
\end{center}
\bigskip
 {\bf Abstract.} We consider selected aspects of (non-Archimedean) quantum mathematics and non-Archimedean (quantum) computation.
\bigskip
\begin{center}
{\bf Motivations and History  of non-Archimedean quantum physics} (very short)
\end{center}
In the summer of 2025, two conferences were held: "On the Resurgence in Quantum Mathematics" at the Centre for Quantum Mathematics, IMADA, University of Southern Denmark and a conference dedicated to the memory of Yu. I. Manin at the Max Planck Institute for Mathematics in Bonn.

Invitations and participation in them were one of the motivations for presenting my report.

{\bf First, the Motivations and History of Arithmetic in Physics (Arithmetic Physics).} (very short) 

We use Number Theory and Arithmetic interchangeably (as synonims).\\
{\t Some history}

From the lecture (ICM 2018) {\bf ARITHMETIC PHYSICS} of Fields Laureate Sir Michael Atiyah:

{\it My Abel lecture in Rio de Janeiro was devoted to solving the riddle of the fine structure
constant $\alpha$, an outstanding problem from the world of fundamental physics.}

 $\alpha \approx 1/137.035999206$ - postoyannaya tonkoj struktury

"In physics, the fine-structure constant, also known as the Sommerfeld constant, commonly denoted by 
 $\alpha$ (the Greek letter alpha), is a fundamental physical constant that quantifies the strength of the electromagnetic interaction between elementary charged particles".
Wiki.

Nobel Laureate  Richard Feynman (about $\alpha$):

{\it Where does $\alpha$ come from; is it related to $\pi$, or perhaps to $e$? Nobody knows, it is one of the great damn mysteries of physics: a magic number that comes to us with no understanding by man. You might say the hand of God wrote the number and we don’t know how He pushed his pencil.}\
\newpage

 Fields Laureate Sir Michael Atiyah:\\
 {\it aims to 
build a Theory of Arithmetic Physics.}

 {\it A programme with a similar title was put forward by Yuri Manin in Bonn in 1984
Atiyah [1985]. In fact Manin was not able to come and I presented his lecture for him.
This gave me the opportunity to suggest that Manin’s programme was too modest, and
that it should be lifted from the classical regime to the quantum regime with Hilbert
spaces in the forefront.}\\

Non-Archimedean (quantum) mathematics began with Hensel's introduction of $p$-adic numbers, which belong to arithmetic.\\

{\it Motivations}

\begin{itemize}
  
\item  {\it  Space-time structure} 

The hypothesis that the geometry of spacetime is non-Archimedean at the sub-Planckian scales.\\

\item  {\it  New models}

The use of non-Archimedean  mathematics provides a new framework for building mathematical models of reality based on different number systems\\

\item  {\it  Unified framework}

Some researchers are seeking a unated framework that combines quantum mathematics and gravity (Quantum Field Theory (QFT), string theory), possibly using non-Archimedean methods

\end{itemize}

\section{$p$-adic numbers} 
$p$-adic numbers introduced at the beginning of the twenties century by Kurt Hensel (1861-1941). \\

  In coding theory and in criptography, codes are constructed and applied over finite fields and rings, for instance, over ${\mathbb F}_p$ and 
over ${\mathbb Z}/p{\mathbb Z}$ for prime $p$  \cite{cs,bgn}.

Under working with codes we have to solve congruences.

Let $n \in {\mathbb N}$. 
\begin{example}
 $x^2 \equiv 2 \mod 7^n$

We have two  solutions $x_0 \equiv3 \mod 7$ and $x_0 \equiv 4 \mod 7$.

Consider the case $x_0 \equiv 3 \mod 7$.

Let $n=2$.

Then $x_1 \equiv 3 + 7\cdot 1 \mod 7^2$.

For $n = 3$ we have: $x_2 \equiv 3 + 7\cdot 1 + 7^2\cdot 2 \mod 7^3$. And so on.
\end{example}

 \begin{definition}
Let $p$ be a prime number. The sequence of integer numbers

 $\{x_0 , x_1, \ldots, x_n, \ldots \}$ 

such that 

$x_n \equiv x_{n-1} \mod p^n$

is called a $p$-adic number.
\end{definition}

\begin{example}
   $216 = 0,0011011 $ (2-adic),\\
    $216 = 0,0022 $ (3-adic),\\
     $216 = 1,331 $ (5-adic)
\end{example}

\section{$p$-adic is not Archimedean}

Let

$a_0 + a_1p + a_2p + \cdots,  \; 0 \le a_i < p.  $

be a $p$-adic integer.

Let ${\mathbb Z}_p = \varprojlim {\mathbb Z}/{p^n{\mathbb Z}}$ be the ring of $p$-adic numbers,\\

\subsection{$p$-adic absolut value}

Let $|\;|_p$ be a $p$-adic absolut value.

For a rational number 
\begin{equation}
\label{pa}
a = p^n \frac{b}{c}, \; (bc,p) = 1
\end{equation}
 let
$|a|_p = (1/p)^n$.

The exponent $n$ in the representation (\ref{pa}) of the number $a$ is denoted $\nu_p(a)$, and one puts $\nu_p(0) = \infty$

Properties of $\nu_p$.
\begin{itemize}
\item $\nu_p(a) = \infty \Longleftrightarrow a=0,$
\item $\nu_p(ab) = \nu_p(a) + \nu_p(b)$,
\item {\it  non-Archimedean property: $\nu_p(a + b) \ge \min \{\nu_p(a), \nu_p(b)\}$} 
\end{itemize}

 The $p$-adic absolut value ($p$-adic norm) is given by
$$
  |\;|_p:  {\mathbb Q} \to {\mathbb R}, \;  a \mapsto \nu_p(a) = 
\left(\frac{1}{p}\right)^{\nu_p(a)},
$$
Properties of $|\;|_p$.
\begin{itemize}
\item $ |a|_p = 0  \Longleftrightarrow a=0,$
\item $|ab|_p = |a |_p |b|_p$,
\item $|a + b|_p \le \max \{ |a|_p,  |b|_p\}$
\end{itemize}
 $|a|_{\infty} = |a|$  (the usual absolute value of a real number)
\subsection{Functional case} Let $k$ be a finite field. We have inside $k(x)$ the polynomial ring $k[x]$ 
      $\mathfrak{p} \not = 0$ of  which are given by the monic irreducible polynomials $p(x) \in k(x)$. For every 
such $\mathfrak{p}$, one defines an absolute value $|\;|_\mathfrak{p}$.

\subsection{Product formula}

For every rational number $a \not = 0$, one has

\begin{equation}
\label{pf}
 \prod_p |a|_p = 1,
\end{equation}
where $p$ values over all prime numbers as well as the symbol $\infty$.
\begin{remark}
A product formula also valid in the functional case.
\end{remark}

\subsection{Adeles}

 \section{$p$-adic codes and their implementation} \cite{kri}
$p$-adic arithmetic with $r$-digits can be considered equivalent to the residue arithmetic modulo $p^r$ \cite{bsh}.\\.

\section{Classical quantum mechanics}
\subsection{Spin in classical quantum mechanics and its properties}

The elementary unit of quantum information is the qubit, a two-state quantum system whose basis states $|0\rangle$ and $|1\rangle$ 
(in Dirac notation) correspond to binary $0$ and $1$.
Quantum information is encoded by arrays of qubits, and elementary transformations of quantum information are controlled by quantum gates.
Group-theoretical methods for transforming quantum information are related to the Pauli and Clifford groups.
The two-dimensional identity matrix $\sigma_0$ and the Pauli matrices $\sigma_0$, $\sigma_x$, $\sigma_y$ $\sigma_z$   form a basis for the space of density operators of one qubit.
For arrays of $n$ qubits, the Pauli group is used
and the Clifford group, defined as the normalizer of the Pauli group.

  \subsection{Spintronics}

\section{Algebraic quantum circuits and exponential sums} 

\subsection{On quantum logic}

Quantum logic is not the same as propositional calculus.
Distributive low is false: $a \wedge (b \vee c) \not = (a \wedge b) \vee (a \wedge c).$

\subsection{Modular structure}
A modular structure is a set $L$  with binary operations $\wedge, \vee,$ that satisfy the following conditions:
\begin{itemize}
\item  {$\wedge, \vee,$  {\it are associative and commutative}}

\item  {$a \wedge a = a \vee a =a$  {\it   for all $a \in L$}}

\item  {{\it If $a \wedge b = b$, then $(a \vee c) \wedge b = a  \vee (c \wedge b)$ (modular identity)  } }

\end{itemize}

\subsection{On quantum circuits} 
Analyzing algebraic quantum circuits using exponential sums (algebraicheskie kvantovye cxemy) \cite{bdr,tsqg,glaap}

\section{Analytification of the field of $p$-adic numbers}

Let ${\mathbb Q}_p$ be the field of $p$-adic numbers,\\

let ${\overline{\mathbb Q}}_p$ be the algebraic closure of ${\mathbb Q}_p$,
 and \\

let ${\mathbb C}_p$ be the complation of the algebraic closure:\\

\begin{center}
${\mathbb Q}_p \to {\overline{\mathbb Q}}_p \to {\mathbb C}_p$
\end{center}

But the field ${\mathbb C}_p$ is not connected. So there are problems with Analysis on it.

Let $A$ be a ring, $K$ be an algebraically closed field which is complete with respect to
a non-archimedean absolute value.

A multiplicative seminorm on a ring
$A$ is a function $|\;|_x: A \to {\mathbb R_{\ge0}}$ satisfying: 
\begin{itemize}
 \item {\it $ |0|_x = 0$ and  $|1|_x = 1$} .

 \item  {\it $ |f g|_x = |f|_x |g|_x$  for all $f, g \in A$}.

 \item {\it $ |f +  g|_x \le  |f|_x + |g|_x$  for all $f, g \in A$}. 

\end{itemize}

\section{Spaces and groups of quantum mathematics}
\subsection{On path integrals}
"It should be noted that the path integral is mathematically strictly defined only in Euclidean space as the integral over the
 Wiener measure." \cite{gg}

"A typical case of Euclidean vector space is 
${\displaystyle \mathbb {R}^{n}}$ viewed as a vector space equipped with the dot product as an inner product. The importance of this particular example of Euclidean space lies in the fact that every Euclidean space is isomorphic to it. More precisely, given a Euclidean space $E$ of dimension $n$, the choice of a point, called an origin and an orthonormal basis of ${\displaystyle {\overrightarrow {E}}}$ defines an isomorphism of Euclidean spaces from $E$ to ${\displaystyle \mathbb {R}^{n}.}$
(Solomentsev, E.D. (2001) [1994], "Euclidean space", Encyclopedia of Mathematics, EMS Press; Wiki.)

\subsection{Quantume groups (V. Drinfeld, Ju. Manin)}

\section{Formal power series methods and non-Archimedean methods in the  theory of resurgence} 
See formal power series methods in the  theory of resurgence in lectures \cite{etr}. 
See a simple example. of a resurgence in the Appendix.

\section{Noncommutative spaces and Alane Connes works}

 \subsection{On A. Connes, C. Consani, H. Moscovici paper
 "Zeta zeros and prolate wave operators : semilocal adelic operators}
The essence of this interesting paper is the integration in the semilocal framework by 
 Connes,  Consani [ J. Oper. Theory 851, 257-276 (2021; Zbl 1524.11171)] 
 of the trace formula by Connes [Selecta Math. (N.S.
 5, 1, 29-106, (1999; Zbl 0945.11015)] two recent results by the first and second authors
[ Enseign. Math., 69, 1-2, 93-148 (2023; Zbl 1537.11121)] and by the first and therd authors  [Proc. Natl. Acad. Sci. USA, 119, 22 (2022; Zbl 07998560) ]
``on the spectral realization of the zeros of the Riemann zeta function by introducing
a semilocal analogue of the prolate wave operator''. 

Based on the (quantum) mathematical physics by Connes, Marcolli [Noncommutative Geometry,
 Quantum Fields and Motives, Colloquium Publications, Providence: American Mathematical
Society, Providence, (2008; Zbl 1209.58007)] the results of the two above-mentioned papers
are interpreted as the infrared (low-lying) part and the ultraviolet behavior of the zeros
 of the Riemann zeta function. 
 
 The main results of the paper under review are Theorem 1 on semilocal Hardy–Titchmarsh transform   and Theorem 2  that relate with ``the sought-for Weil cohomology''.
  The paper  consists of six sections.
   The first section is of introductory nature and designed to acquaint the reader with main ideas and results of the paper.
Section 2 presents cyclic pairs and associated prolate operators.

Let $\mathcal H$ be a Hilbert space, $D$ a selfadjoint  operator acting on 
$\mathcal H$ and $\xi \in \cap_{n \in {\mathbf N}} Dom D^n$ a unit cyclic vector. Let $\lambda > 0$.
For a cyclic pair $(D, \xi)$ authors of the paper under review define the formal prolate operator $\omega(D,\xi,\lambda)$ as 
$\omega(D,\xi,\lambda) = -D^2 + \lambda^2 N$, where $N$  is the  number operator ($N{\xi}_j = j{\xi}_j$).

Let the ${\mathcal C}^*$-algebra 
${\mathcal A} := c_0(\mathbb N)$ be the algebra of sequences 
vanishing at $\infty$.
Authors for an even cyclic pair $(D, \xi)$ and the ${\mathbb Z}/2$ grading $\gamma$ define the even spectral triple $({\mathcal A}, {\mathcal H}, D)$   where the ${\mathcal C}^*$-algebra ${\mathcal A}$ acts in the Hilbert space ${\mathcal H}$.
After a number of preliminary results (Propositions 2.1 - 2.3),
 the authors define the Jacobi matrices of the formal prolate
operator (Proposition 2.4).

Section 3 deals with Hardy–Titchmarsh transform: archimedean place.
 ``Section 4 is devoted to the extension of the Hardy–Titchmarsh transform in the
semilocal situation involving a finite set S of places of Q containing the archimedean
one, and the analysis of the semilocal Sonin space.''
Section 5  studies metaplectic representation in the Jacobi picture.
  The sixth section is the Appendix: Fourier transform on $p$-adic numbers.
  The (Weil) positivity is given in the paper by A. Weil [Comm. Lund., pp. 252–267 (1952; Zbl 0049.03205)]. See also papers by 
 H. Yoshida [In: Zeta Functions in Geometry (Tokyo, 1990), Adv. Stud. Pure Math., vol. 21, pp. 281-325. Kinokuniya, Tokyo (1992; Zbl 0817.11041) ]  and by E Bombieri [Rend. Mat. Acc. Lincei, s. 9, v. 11:183-233 (2000; Zbl )].
 
 The first author recently give the talk at CIRM': Alain Connes, Extremal Eigenvectors, Spectral Action,
and the Zeta Spectral Triple, CIRM, 10 April 2025.

\subsection{On $S$-algebraic group}
E. Lindenstrauss, Invariant measures and arithmetic quantum unique ergodicity.

 by Fields Laureate Lindenstrauss  \cite{lin} "We classify measures on the locally homogeneous space $\Gamma \backslash SL(2, {\bf R})) \times L$
which are invariant and have positive entropy under the diagonal subgroup
of $SL(2, {\bf R}))$ and recurrent under $L$."

by Lindenstrauss and others  "We recall that the group $L$ is $S$-algebraic if it is a finite product of algebraic
groups over ${\bf R}, {\bf C}$, or ${\bf Q}_p$, where $S$ stands for the set of fields that appear in this
product. An $S$-algebraic homogeneous space is the quotient of an $S$-algebraic
group by a compact subgroup."

\subsection{On absolute algebraic geometry}

Alain Connes, Fields Medal in 1982,  operator algebras (von Neumann algebras), absolute algebraic geometry
 and noncommutative geometry.

At first consider an example.
\begin{example}
Consider the Chevalley group $SL_3({\mathbb F}_p)$. This is a simple group of Lie type $A_{2}$.
 Its Weyl group is the symmetric group $S_3$. 
 \end{example}

by \cite{conc}. Categorical object ${\mathbb F}_1$.   A finite field with one element.


\begin{remark}
  Prof. V.A. Ustimenko \cite{ust} explores the problems of post-quantum cryptography.
In his research he develops and applies the results of his studies of symmetric groups and Chevalley groups.
He, in particular, believes (along with other researchers) that for a Chevalley group $G$ one has $G({\mathbb F}_1) = W$, the Weyl group of $G$.
\end{remark}

 Consider the affine scheme $Spec \; {\mathbb Z}$.

Over each  prime finite point $p$ of $Spec \; {\mathbb Z}$ there lies the corresponding prime finite field ${\mathbb F}_p$.

Let $q = p^n, \: n \in {\mathbb N}$. It is possible to consider  ${\mathbb F}_q$ as a quantum deformation of ${\mathbb F}_1$ by $q$.

Over generic point  of $Spec \; {\mathbb Z}$ we have the field of rational numbers  $ {\mathbb Q}$.

Let us construct the projective limit ${\mathbb Z}_p = \varprojlim {\mathbb Z}/{p^n{\mathbb Z}}$ for each  prime finite field ${\mathbb F}_p.$

Let us complete the field of rational numbers $ {\mathbb Q}$ to a real field $ {\mathbb R}$.

On each of the $p$-adic fields, we define the $p$-adic valuations $|\;|_p$ defined above and, on the field of rational numbers, the absolute value. Taking their product, we obtain the left-hand side of formula  (\ref{pf}).


The following problems arise for the field ${\mathbb F}_1$.

\begin{problem}
  Can the field ${\mathbb F}_1$ be completed to a complete field?\\
If so, what is the form of this "$1$-adic complete field"?\\
What is the class of "1-adic norms of this complete field"?
\end{problem}

Recall that there is the Arakelov compactification ${\overline{Spec\:{\mathbb Z}}}$. 

\begin{problem}
If we could modify  the product (\ref{pf}) by adding the object "1-adic norm" to 
the left-hand side of the formula  (\ref{pf}),
 then we could obtain an invariant similar to the product (\ref{pf}). What does it look like?

 The similar for adeles.
 \end{problem} 

\section{Appendix}
\label{rsu}
About resurgence (vozrojdenie):

This term have introduced (defined) by J. Ecalle in 1998 year \cite{etr}.

Consider Euler's equation
\begin{equation}
\label{ee}
 t^2( dy/dt) = t -y
\end{equation}

By a
book  on ODE we try to solve (\ref{ee}) by

$
y = \sum_{m=0}^{\infty} a_m t^m
$

We obtaine

$
 y = \sum_{m=0}^{\infty}(-1)^m m! t^{m+1}
$

Unfortunately this series divergences  for all $t \neq  0$.

But if we consider homogeneous equation 
$  t^2( dy/dt) =  -y  $
and solve it, we have 

$ y = ae^{1/t}.$

The solution
$$
   y = \sum_{m=0}^{\infty}(-1)^m m! t^{m+1} + ae^{1/t}.
$$ 
is the resurgence solution of   (\ref{ee}).

I will call the resurgence the Tauberian resurgence.

There exists also the Abelian resurgence. \\

{\bf Acknowledgments.} I thank Professor Jørgen Ellegård Andersen and Professor Dr. Peter Teichner for inviting me to the conferences.
I thank the leaders and the scientific secretary of the International scientific seminar "Quantum Computing" for the invitation to give a talk.
I would like to thank Academician V. Gusynin for his consultations.

\end{document}